\documentclass{article}
\usepackage{amssymb,latexsym,epsfig}
\begin{document}

\newcommand{\mmbox}[1]{\mbox{${#1}$}}
\newcommand{\proj}[1]{\mmbox{{\bf P}^{#1}}}
\newcommand{\affine}[1]{\mmbox{{\bf A}^{#1}}}
\newcommand{\Ann}[1]{\mmbox{{\rm Ann}({#1})}}
\newcommand{\caps}[3]{\mmbox{{#1}_{#2} \cap \ldots \cap {#1}_{#3}}}
\newcommand{\N}{{\bf N}}
\newcommand{\Z}{{\bf Z}}
\newcommand{\R}{{\bf R}}
\newcommand{\Tor}{\mathop{\rm Tor}\nolimits}
\newcommand{\Ext}{\mathop{\rm Ext}\nolimits}
\newcommand{\Hom}{\mathop{\rm Hom}\nolimits}
\newcommand{\im}{\mathop{\rm Im}\nolimits}
\newcommand{\rank}{\mathop{\rm rank}\nolimits}
\newcommand{\supp}{\mathop{\rm supp}\nolimits}
\newcommand{\arrow}[1]{\stackrel{#1}{\longrightarrow}}
\newcommand{\text}[1]{\mbox{\rm {#1}}}

\newtheorem{defn0}{Definition}[section]
\newtheorem{prop0}[defn0]{Proposition}
\newtheorem{conj0}[defn0]{Conjecture}
\newtheorem{thm0}[defn0]{Theorem}
\newtheorem{lem0}[defn0]{Lemma}
\newtheorem{corollary0}[defn0]{Corollary}
\newtheorem{example0}[defn0]{Example}

\newenvironment{defn}{\begin{defn0}}{\end{defn0}}
\newenvironment{prop}{\begin{prop0}}{\end{prop0}}
\newenvironment{conj}{\begin{conj0}}{\end{conj0}}
\newenvironment{thm}{\begin{thm0}}{\end{thm0}}
\newenvironment{lem}{\begin{lem0}}{\end{lem0}}
\newenvironment{cor}{\begin{corollary0}}{\end{corollary0}}
\newenvironment{exm}{\begin{example0}\rm}{\end{example0}}

\newcommand{\defref}[1]{Definition~\ref{#1}}
\newcommand{\propref}[1]{Proposition~\ref{#1}}
\newcommand{\thmref}[1]{Theorem~\ref{#1}}
\newcommand{\lemref}[1]{Lemma~\ref{#1}}
\newcommand{\corref}[1]{Corollary~\ref{#1}}
\newcommand{\exref}[1]{Example~\ref{#1}}
\newcommand{\secref}[1]{Section~\ref{#1}}

\newcommand{\qed}{\mbox{$\Box$}}
\newenvironment{proof}{\noindent {\bf Proof.}}{\qed\vskip 6pt}
\newenvironment{proofn}{\noindent {\bf Proof.}}{\vskip 6pt}

\newcommand{\poina}{\pi({\cal A}, t)}
\newcommand{\TM}{D({\cal A})}

\parskip = 4pt

\title{Elementary modifications and line configurations in $\mathbb{P}^2$}
\author{Henry K. Schenck\thanks{This work was completed while the author
was an NSF postdoctoral fellow at Harvard University.
\newline   \mbox{   } \mbox{   }   AMS subject classifications 14N20, 14J60, 52C30, 14F05.
\newline   \mbox{   } \mbox{   }   keywords: hyperplane arrangement, vector bundle, 
Castelnuovo-Mumford regularity, stability, jump locus.} \\    
         Texas A\&M University\\
         College Station, TX 77843 \\
         schenck@math.tamu.edu}
\maketitle
\sloppy

\begin{abstract} \noindent 
Associated to a projective arrangement of hyperplanes ${\cal A} \subseteq \mathbb{P}^n$ 
is the module $\TM$, which consists of derivations tangent to ${\cal A}$. We study $\TM$ when 
${\cal A}$ is a configuration of lines in $\mathbb{P}^2$.
In this setting, we relate the deletion/restriction construction used in the study of
hyperplane arrangements to elementary modifications of bundles. This allows
us to obtain bounds on the Castelnuovo-Mumford regularity of $\TM$. We 
also give simple combinatorial conditions for the associated bundle to be stable, 
and describe its jump lines. These regularity bounds and stability considerations
impose constraints on Terao's conjecture.
\end{abstract}


\section{Introduction}\label{sec:intro}

\noindent In this paper we investigate the connection between a
standard construction in algebraic geometry (elementary modifications
of bundles) and a standard construction in the study of hyperplane arrangements
(the deletion/restriction operation). In the setting of line configurations
in $\mathbb{P}^2$, it turns out that they are the same thing. 
Given a two bundle $V$ and line $L$, if $V$ sits in a
modification $$ 0 \longrightarrow W \longrightarrow V \longrightarrow
i_*{\cal O}_L(a) \longrightarrow 0,$$ then understanding $V$ means
understanding $W$ and $a$ (here $i$ is the inclusion of $L$ in $\mathbb{P}^2$). 
For line configurations, the twist $a$ 
has a simple combinatorial meaning, and the long exact sequence
in cohomology yields information about the arrangement. In particular,
it gives a bound on the Castelnuovo-Mumford regularity of the module
of derivations tangent to ${\cal A}$. 

In \S 2 we give a quick review of the fundamental objects:
syzygy modules and hyperplane arrangements. 
The module of ${\cal A}$ derivations is denoted $\TM$, it consists of derivations
of $\mathbb{P}^n$ tangent to ${\cal A}$. The syzygy module on the 
Jacobian ideal of the singular hypersurface ${\cal A}$ is a summand of 
$\TM$, hence the connection. In \S 3 we discuss
elementary modifications and Castelnuovo-Mumford regularity, and
in \S 4 turn our attention to stability. We close with an application of
these results to Terao's conjecture that the freeness of $\TM$
depends only on the combinatorics of the arrangement.

\section{Zero dimensional subschemes and line configurations}\label{sec:hyper}
Let $R = k[x_0, x_1, x_2]$ and $I = \langle s_1,\ldots, s_k \rangle \subseteq
R_m$ a codimension two ideal; $I$ defines a map:
$$ {\cal O}^k_{\mathbb{P}^2} \stackrel{I}{\longrightarrow} {\cal O}_{\mathbb{P}^2}(m).$$
Let ${\cal D}$ denote the sheaf associated to the module of syzygies
on $I$. ${\cal D} \simeq \oplus{\cal O}(\gamma_i)$ iff $I$ is saturated. Since
${\cal D}$ is a second syzygy sheaf, ${\cal D}$ is locally free. Let
$Z$ be the scheme defined by $I$; the following lemma follows from
standard properties of Chern classes:
\begin{lem}\label{lem:cclass}
The Chern classes of ${\cal D}$ are $c_1({\cal D})=-m$, $c_2({\cal D}) = m^2-deg(Z)$.
\end{lem}

\noindent In \cite{se}, Serre describes a method of
constructing a rank two vector bundle ${\cal F}$ on $\mathbb{P}^n$ from a
codimension two local complete
intersection $Y$ with ideal sheaf ${\cal I}_Y$. If the determinant
bundle of the normal bundle of $Y$ extends to a bundle on $\mathbb{P}^n$:
$$\mbox{det }{\cal N}_{Y/\mathbb{P}^n} \simeq {\cal O}_{\mathbb{P}^n}(m)|_Y,$$
then there is a rank two bundle ${\cal F}$ on $\mathbb{P}^n$ with section $s$,
which induces the short exact sequence:
$$0 \longrightarrow {\cal O}_{\mathbb{P}^n} \stackrel{\cdot s}{\longrightarrow}
{\cal F} \longrightarrow {\cal I}_Y(m) \longrightarrow 0.$$ 
The Chern classes of ${\cal F}$ are given by $c_1({\cal F})=m$ 
and $c_2({\cal F})= \mbox{deg }Y$. If $Z$ is a local complete intersection and the bundle
${\cal F}$ exists, then the bundles ${\cal D}$ and
${\cal F}$ are related by the exact sequence:
$$0 \longrightarrow {\cal D} \longrightarrow {\cal O}^{k+1} \longrightarrow
{\cal F} \longrightarrow 0.$$
Let char $k = 0$ and let $Q\in R_{m+1}$ be a reduced polynomial; the role of
$I$ will be played by the Jacobian ideal of $Q$. For the remainder of the 
paper we restrict our attention to the case where $Q$ is a product of
distinct linear forms, although many of the results can be generalized. 
We begin with some facts about hyperplane arrangements; for more
information  see Orlik and Terao \cite{ot}.

\noindent A hyperplane arrangement
${\mathcal A}$ is a finite collection of codimension one linear subspaces
of a fixed vector space V.
${\mathcal A}$ is {\it central} if each hyperplane contains the origin
{\bf 0} of V.
The intersection lattice $L_{\mathcal A}$ of ${\mathcal A}$ consists of the
intersections of
the elements of ${\mathcal A}$; the rank of $x \in L_{\mathcal A}$ is 
simply the codimension of $x$. V is the
lattice element $\hat{0}$; the rank one elements are the hyperplanes
themselves. ${\mathcal A}$ is called {\it essential} if rank $L_{\mathcal
A} =$ dim $V$. Henceforth,
${\cal A}$ will be an {\it essential}, {\it central} three arrangement with $|{\cal A}|=d$; i.e.
a set of $d$ lines in $\mathbb{P}^2_k$.

\begin{defn}
The M\"{o}bius function $\mu$ : $L_{\mathcal A} \longrightarrow \mathbb{Z}$ is
defined
by $$\begin{array}{*{3}c}
\mu(\hat{0}) & = & 1\\
\mu(t) & = & -\sum\limits_{s < t}\mu(s) \mbox{, if } \hat{0}< t
\end{array}$$
\end{defn}

\begin{defn}
The Poincar\'e polynomial $\poina = \sum\limits_{x \in L_{\mathcal
A}}\mu(x) \cdot (-t)^{rank(x)}.$
\end{defn}

\noindent  It follows directly
from the definitions above that for such an arrangement, $$\poina =
(1+t)(1+(d-1)t+(\sum\limits_{\stackrel{x \in L_{\cal A}}{rank(x)=2}}\mu(x) -d+1)t^2).$$
Let $Q$ be a reduced polynomial defining ${\cal A}$ and $J_Q$ the Jacobian
ideal of $Q$. The next lemma gives an easy proof of the main result of \cite{sch}:
\begin{lem}\label{lem:lci}
The Jacobian ideal of a line arrangement in $\mathbb{P}^2$ is a local complete intersection.
\end{lem}
\begin{proof}
Localization, the product rule, and Euler's relation.
\end{proof}
\begin{thm}\label{thm:compute1} If ${\mathcal A}$ is a line
arrangement defined by $Q$ and ${\cal D}$ is the syzygy bundle of 
$J_Q$, then
$$\poina = (1+t) \cdot c_t({\cal D}^\vee),$$
where $c_t$ is the Chern polynomial.
\end{thm}
\begin{proof}
\lemref{lem:cclass} implies that $$c_t({\cal D}^\vee) = 1 + (d-1)t +
((d-1)^2 - \mbox{deg }J_Q)t^2,$$
and  by \lemref{lem:lci} we have  $$\mbox{deg }J_Q = \sum\limits_{\stackrel{x
\in L_{\cal A}}{rank(x)=2}}\mu(x)^2.$$ 
Now use the identity: $${d \choose 2} = \sum\limits_{\stackrel{x \in L_{\cal
A}}{rank(x)=2}}{\mu(x)+1 \choose 2}.$$
\end{proof}
The motivation for the previous theorem is Terao's celebrated freeness
theorem \cite{t}. Let ${\cal A} \subseteq \mathbb{P}^n$ be an arrangement with
defining polynomial $Q \in R = k[x_0, \ldots, x_n]$;  the module of derivations
tangent to ${\cal A}$ is defined as:
\begin{defn} $\TM = \{ \theta \in Der_k(R)\mbox{ }  | \mbox{ } \theta(Q) \in
\langle Q \rangle\}.$
\end{defn}
Terao's theorem is that if $D({\cal A})$ is free, then the Poincar\'e polynomial  
factors as $\Pi_{i=0}^n(1+a_it)$, where $a_i$ are the degrees of a set
of homogeneous generators of $D({\cal A})$. 
If char $k = 0$ then $\TM \simeq D_0 \oplus R(-1)$,
where $R(-1)$ is generated by the Euler derivation and
$D_0$ is the module of syzygies on $J_Q$. Henceforth, ${\cal D}$
will be the sheaf associated to $D_0$.

\section{Elementary Modifications and Castelnuovo-Mumford regularity}\label{sec:elm}
\noindent A triple of arrangements $({\cal A}', {\cal A}, {\cal A}'')$
consists of an arrangement ${\cal A}$
and choice of distinguished hyperplane $H \in {\cal A}$ such that ${\cal
A}' = {\cal A} - H$ and ${\cal A}'' = {\cal A}|_H$.
${\cal A}' $ is called the deletion of ${\cal A}$ and ${\cal A}''$ the restriction
of ${\cal A}$ with respect to $H$. Of course, the invariants of the elements
of a triple are closely related, and deletion-restriction is often a valuable
tool for inductive proofs. For the module of derivations, we have:
\begin{prop}\label{prop:ot}$($Proposition 4.45 of \cite{ot}$)$  There is an
exact sequence:
$$ 0 \longrightarrow D({\cal A}')(-1) \stackrel {\cdot H}{\longrightarrow}
D({\cal A}) \longrightarrow
D({\cal A}'')$$
\end{prop}
Orlik and Terao give an example of a line arrangement for which the above
sequence is not right exact (see example 4.56 of \cite{ot} - this
corresponds to example I in the next section, with the role
of $H$ played by $\{z=0\}$). Recall the definition of an elementary
modification (see
\cite{fr}): Let $X$ be a smooth variety, $Y$ an effective divisor on
$X$,
$Y \stackrel{i}{\hookrightarrow} X$. Let $V$ be a rank two bundle on $X$,
$L$ a line bundle
on $Y$, and suppose $V \rightarrow i_*L \rightarrow 0$. Then the kernel $W$
of the map is
also a rank two bundle on $X$, with $c_1(W) = c_1(V)-Y$ and $c_2(W) =
c_2(V)-c_1(V)\cdot Y+i_*c_1(L)$.
\begin{thm}\label{thm:elem} Let ${\cal A}$ be an arrangement of lines in
$\mathbb{P}^2$. If $({\cal A}', {\cal
A}, {\cal A}'')$ is a triple, and $i: H \simeq \mathbb{P}^1 \hookrightarrow
\mathbb{P}^2$, then the sequence
of sheaves corresponding to \propref{prop:ot} $($with Euler derivations
pruned off$)$ is an elementary modification, i.e.
$$ 0 \longrightarrow {\cal D}'(-1) \longrightarrow {\cal D} \longrightarrow
i_*{\cal D}'' \longrightarrow 0$$
is exact.
\end{thm}
\begin{proof}
This follows since $$i_*{\cal D}''\simeq {\cal O}_H(1-|{\cal A}''|). $$
Now use the Hirzebruch-Riemann-Roch theorem and Serre vanishing to convert the
Chern polynomials to Hilbert polynomials, and compute.
\end{proof}
An important measure of the complexity of a coherent sheaf on $\mathbb{P}^n$ is
the {\it Castelnuovo-Mumford regularity}:
\begin{defn}
A coherent sheaf ${\cal F}$ on $\mathbb{P}^n$ is $m$-regular $(reg({\cal
F})=m)$ if
$$H^i({\cal F}(m-i)) = 0 \mbox{ }\forall  i\ge 1.$$
\end{defn}
The exact sequence $ 0 \longrightarrow {\cal D}'(-1) \longrightarrow
{\cal D} \longrightarrow i_*{\cal D}'' \longrightarrow 0$ gives us a
good way to bound the Castelnuovo-Mumford regularity of line arrangements. 
\begin{thm}\label{thm:cmreg}
For a triple $({\cal A}', {\cal A}, {\cal A}'')$ of line arrangements,
$$reg({\cal D}) \le max \{reg({\cal D}')+1, |{\cal A}''|-1\}.$$
\end{thm}
\begin{proof}
Follows from the long exact sequence in cohomology.
\end{proof}
\begin{cor}
For an arrangement on $d$ lines,  $reg({\cal D}) \le d-2$.
This is tight for generic arrangements.
\end{cor}
Since the sheaf ${\cal D}$ is reflexive, it corresponds to a bundle on $\mathbb{P}^2$ and (\cite{ms})
$$D_0 \simeq \bigoplus\limits_{i}H^0({\cal D}(i)).$$
For an arrangement of $d$ lines, this means that the minimal free resolution
of $D_0$ is:
$$ 0 \longrightarrow \bigoplus\limits_{j=1}^{m-2} R(-\beta_j) \longrightarrow 
\bigoplus\limits_{i=1}^{m}  R(-\alpha_i) \longrightarrow D_0 \longrightarrow 0,$$
where the $\alpha_i$ are at most $d-2$ and the $\beta_j$ are at most $d-1$.
In \cite{z1}, Ziegler gives bounds on the degrees of generators for $D_0^\vee$,
which gives bounds on the generators of $D_0$. For hypersurfaces 
with only isolated singularities, Choudary and Dimca \cite{cd} give 
a bound; for a (reduced, singular) degree $d$ curve in $\mathbb{P}^2$, 
the regularity of $D_0$ is at most $2d-4$. Thus, for line arrangements, 
the bound above is better than existing results.

The entire free resolution of
$D_0$ for a {\it generic} arrangement is given in \cite{rt} and \cite{y}. In \cite{y3},
Yuzvinsky gives a set of generators for a submodule of $D_0^\vee$; these
generators are determined by $L_2({\cal A})$. When the third relation space
vanishes, they actually generate the entire module (this generalizes the results
of Ziegler mentioned earlier). Even in the case of line configurations,
there are examples (e.g. the Braid arrangement) where this space does not vanish.
However, the maximal number of generators of $D_0$ is bounded by $d-1$
(which is attained by generic arrangements).
For one proof of this, see Jiang and Feng \cite{jf}, \S 4.2. 
Finally, we note that Derksen and Sidman \cite{ds} have recently 
obtained regularity bounds on $\TM$ for higher dimensional arrangements. 

\section{Stability and Jump Loci}\label{sec:stab}
In this section, we consider the stability and jump loci of the bundle ${\cal D}$
obtained from an arrangement of $d$ lines; the point is that we can often construct
stable bundles with prescribed jumping lines. For generic arrangements these
questions were studied by Dolgachev and Kapranov in \cite{dk2},\cite{dk1}.
We want to investigate ${\cal D}$ when the arrangement is nongeneric; 
the tool will be the short exact sequence of the last section. We
first recall a few standard results about vector bundles on $\mathbb{P}^n$,
referring
for proofs to the book of Okonek, Schneider and Spindler \cite{oss}.
\begin{defn}
Let ${\cal M}$ be a bundle on $\mathbb{P}^n$. The slope of ${\cal M}$ is defined
as 
$$slope(M) =  \frac{c_1({\cal M})}{\mbox{rk }{\cal M}}.$$
\end{defn}
A key concept in the study of bundles on $\mathbb{P}^n$ is {\it stability}:
\begin{defn}
A bundle ${\cal M}$ on $\mathbb{P}^n$ is stable  
if for all subsheaves ${\cal N} \subseteq {\cal M}$ with $0 < {\mbox{rk }{\cal N}} <{\mbox{rk }{\cal M}}$,
$$slope({\cal N})  <  slope({\cal M}),$$
and semistable if 
$$slope({\cal N})  \le slope({\cal M}).$$
\end{defn}
If ${\cal M}$ is a stable two bundle on $\mathbb{P}^2$, then Schwarzenberger \cite{schw}
showed that $c_1({\cal M})^2 < 4c_2({\cal M})$ (this was generalized
by Bogomolov \cite{bog}). Thus,  a necessary condition for stability of ${\cal D}$ is that
degree $J_Q < \frac{3}{4}(d-1)^2$.
The degree of the Jacobian ideal ranges from ${d \choose 2}$ (achieved for
generic line arrangements)
to $d^2-3d+3$ (achieved for arrangements with $d-1$ lines through a point, and
one other line
in general position, this class of arrangements has ${\cal D} \simeq {\cal
O}(-1) \oplus {\cal O}(-d+2)$).
Basically, as the degree of the Jacobian ideal gets large, ${\cal D}$ has
less chance of being
stable.

\noindent For generic arrangements ${\cal A}$ in $\mathbb{P}^n$, Dolgachev and
Kapranov  prove that the
bundle of meromorphic one forms with logarithmic pole along ${\cal A}$ is
stable (in the setting of generic line arrangements this bundle is dual 
to a twist of ${\cal D}$, see \cite{ms}), they
 also prove that the map which associates 
to a set of $d$ generic hyperplanes the corresponding
bundle is generically injective if $d \ge 2n+3$. Since dim ${\cal M}_{\mathbb{P}^2}(c_1,c_2) = 4c_2-c_1^2-3$
this map has no chance of having Zariski dense image if the number of lines is
large (if $d = 6$ Dolgachev and Kapranov show that it is). First, a few
definitions. The normalization ${\cal F}_{norm}$ of a rank two bundle
${\cal F}$  is ${\cal F}(i)$, where $i$ is chosen so that $c_1({\cal F}(i))
\in \{0,-1\}$, i.e.
\begin{defn}
Let ${\cal F}$ be a rank two bundle on $\mathbb{P}^2$. Then ${\cal F}_{norm} =
{\cal F}(k_{\cal F})$, where $k_{\cal F} =
\frac{-c_1({\cal F})}{2}$ if $c_1({\cal F})$ is even, and $-\frac{c_1({\cal
F})+1}{2}$ if $c_1({\cal F})$ is odd.
\end{defn}
\begin{lem}\label{lem:norm}
A reflexive sheaf ${\cal F}$ of rank two over ${\mathbb{P}^n}$ is stable iff
$H^0({\cal F}_{norm}) = 0$.
If $c_1({\cal F})$ is even, then ${\cal F}$ is semistable iff $H^0({\cal
F}_{norm}(-1)) = 0$. If $c_1({\cal F})$ is odd, semistable and stable coincide.
\end{lem}

\begin{thm}\label{thm:stab} Let ${\cal A}$ be an arrangement of $d$ lines
in $\mathbb{P}^2$, $H$ any
line in ${\cal A}$, and $({\cal A}', {\cal A}, {\cal A}'')$ the
corresponding triple $($notice that the restriction to $H$ ignores 
multiplicities, so $|{\cal A}''|= |{\cal A}' \cap H|)$. Then
\begin{enumerate}
\item If $d$ is odd, then ${\cal D}$ is stable if ${\cal D}'$ is stable and $|{\cal A}''| > \frac{d+1}{2}$.

\item If $d$ is odd, then ${\cal D}$ is semistable if ${\cal D}'$ is semistable and $|{\cal A}''| > \frac{d-1}{2}$.

\item If $d$ is even, then ${\cal D}$ is stable if ${\cal D}'$ is semistable and $|{\cal A}''| > \frac{d}{2}$.
\end{enumerate}
\end{thm}
\begin{proof}
The proofs are similar so we prove the first statement. Since $d$ is odd,
$k_{\cal D}=\frac{d-1}{2}$
and $k_{{\cal D}'}=k_{\cal D}-1$. By \lemref{lem:norm} ${\cal D}$ is stable
iff $H^0({\cal D}_{norm}) =0$ so
it is necessary that $H^0({\cal D}'(-1)\otimes{\cal O}(k_{\cal D})) =
H^0({\cal D}'_{norm}) = 0$. The
simplest sufficiency criterion is then that $H^0(i_*{\cal D}''\otimes{\cal
O}(k_{\cal D})) = 0$; from
the exact sequence
$$0 \longrightarrow {\cal O}(-|{\cal A}''|) \longrightarrow {\cal
O}(1-|{\cal A}''|) \longrightarrow
i_*{\cal D}''
\longrightarrow 0$$ we have that $H^0(i_*{\cal D}''\otimes{\cal O}(k_{\cal D})) = 0$
 iff $|{\cal A}''| > \frac{d+1}{2}$.
\end{proof}
An iff criterion for stability can be formulated
in the obvious fashion: ${\cal D}$ is stable iff $H^0({\cal D}'(-1)\otimes{\cal O}(k_{\cal D})) = 0$
and the connecting map $H^0(i_*{\cal D}''\otimes{\cal O}(k_{\cal D}))
\longrightarrow H^1({\cal D}'(-1)\otimes{\cal O}(k_{\cal D}))$ is an inclusion.
This can be computed in any particular instance, but is not a criterion 
which is easy to apply; the point of \thmref{thm:stab} is that if one
adds a line $H$ to a configuration ${\cal A}'$ for which ${\cal D}'$ is
stable, then as long as no more than (roughly) half of the lines 
become redundant on restricting to $H$, the new configuration is also stable.

Next, we study the jump loci of these bundles. 
By a theorem of Grothendieck, any bundle on $\mathbb{P}^1$ splits as a direct
sum of line bundles;
in particular, for a rank two bundle $V$ on $\mathbb{P}^2$ and line $L$
$$V|_L \simeq {\cal O}_L(a_1)\oplus {\cal O}_L(a_2).$$
In the Grassmannian of lines on $\mathbb{P}^2 \simeq \mathbb{P}^{2^\vee}$, there
is a nonempty, Zariski
open subset where $(a_1,a_2)$ is constant; the complement of this subset is
the jump locus $j_V$ of $V$. By semicontinuity and the
Grauert-M\"ulich
theorem, we have
\begin{defn} For a normalized, semistable two bundle $V$ on $\mathbb{P}^2$
$$j_V = \{ L \in \mathbb{P}^{2^\vee}\mbox{ } | \mbox{ }H^0(V(-1)|_L) \ne 0 \}.$$
\end{defn}
We now describe the jump lines of ${\cal D}$. The method works in
general if we have a modification $V \longrightarrow i_*{\cal O}_L(a)
\longrightarrow
0$. In, \cite{bar}, Barth proved that if $V$ is a semistable two bundle 
on $\mathbb{P}^n$ with $c_1(V) = 0$, then $j_V$ is purely one
codimensional, and is supported on a divisor of degree $c_2(V)$.
If $c_1(V)=-1$, then in general $V$ has only a finite number of jump lines;
in  \cite{hul} Hulek introduced the notion of a jump line of the second kind:
$$j^2_V = \{ L \in \mathbb{P}^{2^\vee}\mbox{ } | \mbox{ }H^0(V|_{L^2}) \ne 0 \},$$
and proved that if $V$ is a semistable two bundle on $\mathbb{P}^n$ 
with $c_1(V) = -1$, then $j^2_V$ is a
curve of degree  $2(c_2(V)-1)$, and $j_V \subseteq \mbox{Sing}(j^2_V).$
We begin by asking which of the lines of the arrangement are contained
in $j_{\cal D}$. Let $({\cal A}', {\cal A},{\cal A}'')$ be a triple with 
respect to a line $H$, put $L = H$.
\begin{thm}\label{thm:jloci1}Let ${\cal A}$ be an arrangement of $d$
lines, $L \in {\cal A}$ and ${\cal A}'' = {\cal A}|_L$.

If $d$ is odd, $L \in j_{\cal D}$ iff either $|{\cal A}''| \ge \frac{d+3}{2}$ or
det $\alpha_L = 0$;
where $$H^0({\cal O}_L\left(\frac{d-1}{2}- |{\cal A}''|\right))
\stackrel{\alpha_L}{\rightarrow}
H^1({\cal O}_L\left(|{\cal A}''| - \frac{d+3}{2}\right)).$$

If $d$ is even, $L \in j_{\cal D}$ iff either $|{\cal A}''| \ge \frac{d+4}{2}$
or rank $\alpha_L <\frac{d}{2}
-|{\cal A}''|$; where $$H^0({\cal O}_L\left(\frac{d-2}{2}- |{\cal A}''|\right))
\stackrel{\alpha_L}{\rightarrow}
H^1({\cal O}_L\left(|{\cal A}''| - \frac{d+4}{2}\right)).$$
\end{thm}
\begin{proof}
We prove the first statement. Since $$i_*{\cal D}''(k_{\cal D} -1)|_L
\simeq {\cal D}''(k_{\cal D}-1)
\simeq {\cal O}_L\left(\frac{d-1}{2} - |{\cal A}''|\right),$$
and $a_1 + a_2 = -2$, we have an exact sequence:
$$0 \longrightarrow {\cal O}_L\left(|{\cal A}''| -
\frac{d+3}{2}\right)\longrightarrow {\cal O}_L(a_1) \oplus  {\cal
O}_L(a_2) \longrightarrow {\cal O}_L\left(\frac{d-1}{2} - |{\cal A}''|\right) \longrightarrow 0.$$ From the long exact sequence in cohomology, it is
obvious that $H^0({\cal D}_{norm}(-1)|_L) \ne 0$ iff one of the two conditions of the theorem holds.
\end{proof}
\pagebreak
\begin{exm} Consider the following set of arrangements (where $\{z=0\}$ is
the line at infinity):

\begin{center}
\epsfig{file=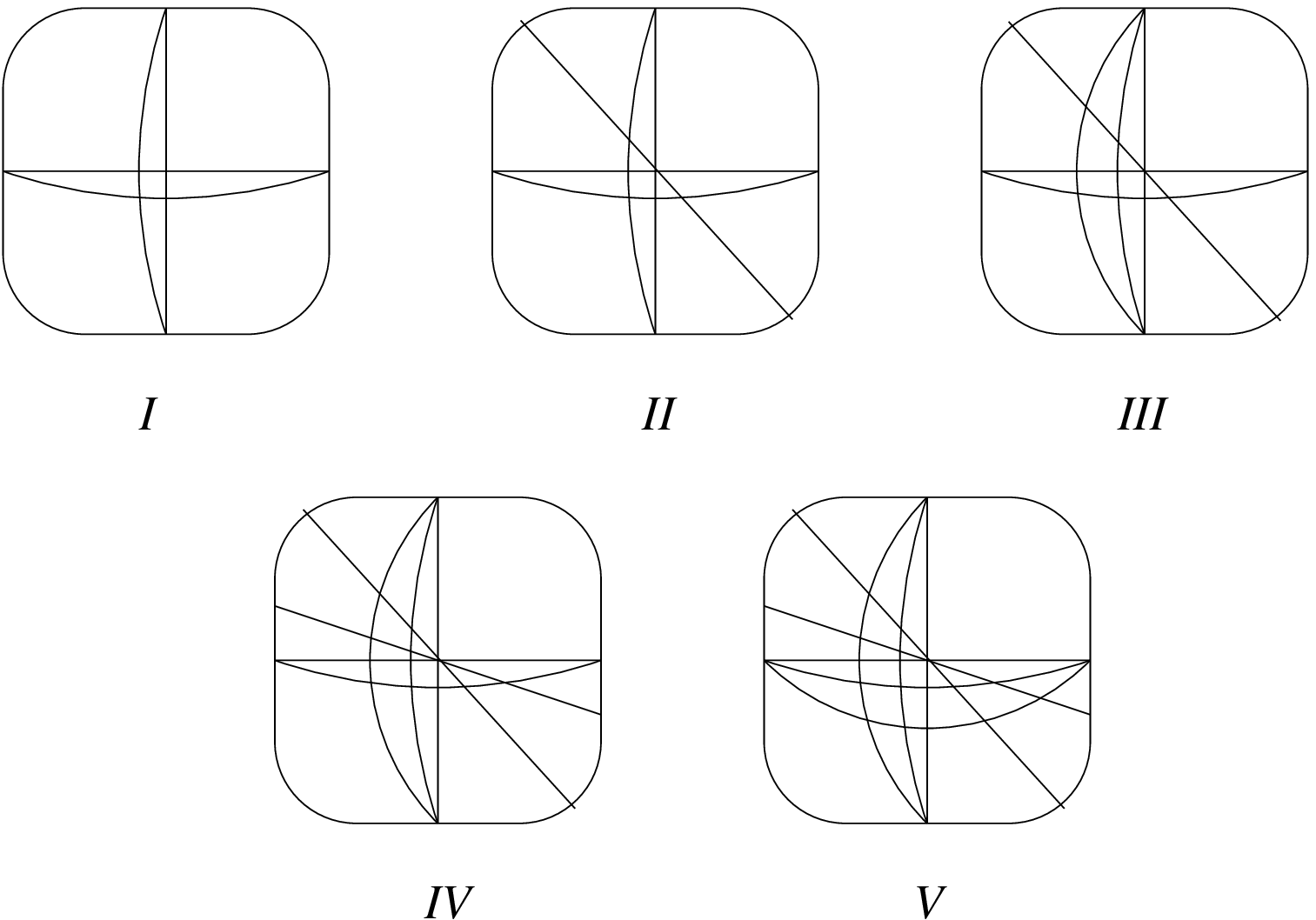,height=3.0in,width=4.3in}
\end{center}

\noindent Arrangement I consists of the coordinate lines and the lines
$\{y+z=0 \}, \{x+z=0 \}$;
${\cal D} \simeq {\cal O}(-2)^2$ (which can be proved using Theorem 4.51 of
\cite{ot}) so ${\cal D}$ is
semistable.
Arrangements II through V are obtained by adding (successively) the lines
$\{x+y=0 \},
\{x+2z=0 \},\{x+2y=0 \},\{y+2z=0 \}$. \thmref{thm:stab} implies that the
bundles associated to
arrangements II, IV, and V are stable and the bundle associated to
arrangement III is semistable. We can compute the Chern classes 
using \thmref{thm:compute1}:
$$\begin{array}{*{5}c}
 & &c_1({\cal D}_{norm}) & &c_2({\cal D}_{norm})\\
I. & & 0 & &0\\

II. & &-1 & &1\\

III. & & 0 & &1\\

IV. &  & -1 & & 2\\

V. &  & 0 & &3
\end{array}$$
The normalized bundles fit into exact sequences:
$$\begin{array}{*{10}c}
I. &  &  & 0 & \longrightarrow & {\cal O}^2 &
\longrightarrow & {\cal D}_{norm} & \longrightarrow & 0 \\

II. & 0 & \longrightarrow & {\cal O}(-2)& \longrightarrow & {\cal O}(-1)^3 &
\longrightarrow & {\cal D}_{norm} & \longrightarrow & 0 \\

III. & 0 & \longrightarrow & {\cal O}(-2)& \longrightarrow & {\cal O}(-1)^2
\oplus{\cal O}&
\longrightarrow & {\cal D}_{norm} & \longrightarrow & 0 \\

IV. & 0 & \longrightarrow & {\cal O}(-3)& \longrightarrow & {\cal O}(-1)^2
\oplus{\cal O}(-2)&
\longrightarrow & {\cal D}_{norm} & \longrightarrow & 0 \\

V. & 0 & \longrightarrow & {\cal O}(-3)& \longrightarrow & {\cal O}(-1)^3 &
\longrightarrow & {\cal D}_{norm} & \longrightarrow & 0 
\end{array}$$
\noindent An easy application of the Beilinson spectral sequence \cite{oss} shows that
a stable bundle with $c_1 = -1$ and $c_2 = 1$ must be $\Omega^1_{\mathbb{P}^2}(1)$, so for
arrangement II, the jump
locus of ${\cal D}$ is empty. For arrangement III, Barth's theorem implies
that the support of $j_{\cal D}$
is of degree one; by \thmref{thm:jloci1}
$\{x+y=0 \}$ and $\{y+z=0\}$ are jumping lines, hence $j_{\cal D}$ is the
line $\{x-y+z=0\} \subseteq \mathbb{P}^{2^\vee}$. For arrangement IV, \thmref{thm:jloci1} shows that
$\{y+z=0\}$ is a jumping line. In fact
(Hulek, Prop 8.2) a semistable two bundle $V$ with $c_1(V) = -1$, $c_2(V) =
2$ has only a single
jump line, so this is it!

\noindent Finally, for arrangement V, Barth's theorem implies that the jump
locus
is a cubic curve; \thmref{thm:jloci1} gives six lines of the arrangement
which are jump lines.
A computation shows that $\{x-y=0 \}$, $\{x-2z=0\}$, and $\{y-z=0\}$ are
also jump lines, so we have
nine points in $\mathbb{P}^{2^\vee}$, which unfortunately only impose eight
conditions on cubics. However,
a final computation shows that the coordinate lines are not jump lines,
which allows us to determine
that the jump locus is a smooth cubic curve:
$$4x^3-2x^2y-4xy^2+2y^3-4x^2z -y^2z-xz^2-2yz^2+z^3 = 0.$$
\end{exm}
\noindent In \cite{dk1}, Dolgachev and Kapranov show that for generic
arrangements in any
dimension, a line contained in one of the hyperplanes of the arrangement is
a jumping line;
for line arrangements this also follows from \thmref{thm:jloci1}. In fact,
Dolgachev and
Kapranov prove that for a generic arrangement with an odd number of lines,
the points of
the jump locus corresponding to the lines of the arrangement are singular
points of fairly high
multiplicity, which is an interesting contrast to the example above.

If $L \ne H$, then since
${\cal T}or_1(i_*{\cal D}'', {\cal O}_L)=0$, we obtain an exact sequence:
$$0 \longrightarrow {\cal D}'(-1)|_L \longrightarrow {\cal D}|_L\longrightarrow i_*{\cal D}''|_L \longrightarrow 0,$$
with $i_*{\cal D}''|_L$ torsion. If $d$ is odd, normalizing, restricting to
$L$, and taking cohomology shows that $j_{{\cal D}'} \subseteq j_{{\cal D}}.$ 
In fact, we can do better.
Let $L$ be a line which is not a line of the arrangement. Think of the
arrangement as ${\cal A}'$, $|{\cal A}'|=d$, and put ${\cal A} = {\cal A}' \cup L$. Restricting to $L$ yields a long exact sequence
$$0 \longrightarrow {\cal T}or_1(i_*{\cal D}'', {\cal O}_L) \longrightarrow
{\cal D}'(-1)|_L \longrightarrow {\cal D}|_L\longrightarrow i_*{\cal D}''|_L \longrightarrow 0.$$
In the previous theorem we studied a short exact sequence pruned from 
the right end of this type of complex; pruning a short exact sequence 
from the left end yields exact sequences:\newline
$$0 \longrightarrow {\cal O}_L\left(\frac{d-1}{2} - |{\cal A}''|\right)\longrightarrow {\cal O}_L(b_1) \oplus  {\cal O}_L(b_2) 
\longrightarrow {\cal O}_L\left(|{\cal A}''|-\frac{d+3}{2}\right) 
\longrightarrow 0,$$
when $d$ is odd; in this case $b_1+b_2=-2$. 
$$0 \longrightarrow {\cal O}_L\left(\frac{d-2}{2} - |{\cal A}''|\right)\longrightarrow {\cal O}_L(a_1) \oplus  {\cal O}_L(a_2) 
\longrightarrow {\cal O}_L\left(|{\cal A}''|-\frac{d+4}{2}\right) 
\longrightarrow 0,$$ 
when $d$ is even; in this case $a_1+a_2=-3$. \newline
\begin{thm}\label{thm:jloci2}
Let ${\cal A}$ be an arrangement of $d$ lines, $L$ a line with $L
\notin {\cal A}$ and ${\cal A}'' = {\cal A}|_L$.\newline

If $d$ is odd, then $L \in j_{\cal D}$ iff either $\frac{d-1}{2} \ge |{\cal A}''|$ or det $\alpha_L = 0$;
where $$H^0({\cal O}_L\left(|{\cal A}''|-\frac{d+3}{2}\right))
\stackrel{\alpha_L}{\rightarrow}  H^1({\cal O}_L\left(\frac{d-1}{2}- |{\cal
A}''|\right)).$$

If $d$ is even, then $L \in j_{\cal D}$ iff either $\frac{d-2}{2} \ge |{\cal A}''|$ or 
rank $\alpha_L \le (|{\cal A}''|-\frac{d+4}{2})$;
where $$H^0({\cal O}_L\left(|{\cal A}''|-\frac{d+4}{2}\right))
\stackrel{\alpha_L}{\rightarrow}  H^1({\cal O}_L\left(\frac{d-2}{2}- |{\cal A}''|\right)).$$
\end{thm}
Notice that the three jump lines of arrangement $V$ which are not lines of 
the arrangement are characterized by the first condition of the theorem. 
It seems reasonable to expect that the jump locus of the second kind is
related to {\it multiarrangements}, about which very little is known
(see Ziegler \cite{z}, or Solomon-Terao \cite{st} for recent progress). 
We plan to return to this question in a later paper.

\section{Terao's conjecture}\label{sec:tconj}
One of the major open conjectures in the study of hyperplane arrangements
is the following:
\begin{conj}$($Terao$)$
In characteristic zero, freeness of $\TM$ depends only on the combinatorics of 
${\cal A}$.
\end{conj}
In \cite{y1}, Yuzvinsky proves that for a fixed intersection
lattice the set of free arrangements is Zariski open. In this section, 
we show that the vector bundle viewpoint has
implications for Terao's conjecture. As noted earlier, $\TM$ decomposes 
as $R(-1) \oplus D_0$, and the module
$D_0$ has a minimal free resolution of the form (with $m \le d-1$):
$$ 0 \longrightarrow \bigoplus\limits_{j=1}^{m-2} R(-\beta_j) \longrightarrow 
\bigoplus\limits_{i=1}^{m}  R(-\alpha_i) \longrightarrow D_0 \longrightarrow 0.$$
If there is an $\alpha_i = 1$, then the arrangement is a ``near pencil'', and
we ignore this case, since such an arrangement is supersolvable \cite{ot}. The results of \S 3 imply
$$2 \le \alpha_i \le d-2 $$
$$3 \le \beta_j \le d-1.$$
Now fix an arrangement with intersection lattice $L_{\cal A}$ (which also fixes
$\pi({\cal A},t)$). For Terao's conjecture we'll be interested in the situation 
where $$\pi({\cal A},t)=(1+t)(1+at)(1+(d-1-a)t).$$ The results on the
Jacobian ideal in \S 2 impose the additional conditions on the resolution:
$$\sum \alpha_i - \sum \beta_j = d-1$$
$$\mbox{deg }J_Q = \sum\limits_{\stackrel{x \in L_{\cal
A}}{rank(x)=2}}\mu(x)^2 = (d-1)^2-a(d-1-a).$$
An easy computation shows the condition on $\mbox{deg }J_Q$ is
equivalent to
$${d-2 \choose 2} -\sum\limits_{i=1}^m {\alpha_i-1 \choose 2}
+\sum\limits_{j=1}^{m-2} {\beta_j-1 \choose 2} +1=a(d-1-a).$$
There are constraints on any free resolution of
the form above; for example, if $\alpha_1 = \min\{ \alpha_i \}$ then
$$\min\{ \beta_i \} \ge \min \{ \alpha_2, \ldots, \alpha_m \}+1.$$
This follows since any relation involves at least two generators, and must
be a positive degree multiple of both. So if a counterexample to Terao's
conjecture exists, it must be an integral solution of the above inequalities.
Another constraint is the following, which follows easily by localization.
\begin{lem}\label{lem:bigmu}
There is a syzygy of degree $\ge \max \{ \mu(x)  |  x \in L_2({\cal A}) \}$.
\end{lem}
\begin{exm}
Let ${\cal A}$ be an arrangement consisting of five lines through a point, 
and two additional lines in general position. 

\begin{center}
\epsfig{file=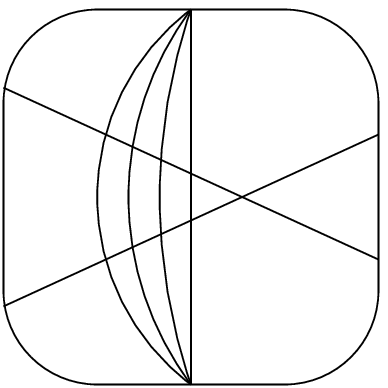,height=1.1in,width=1.1in}
\end{center}
A quick check shows that
$\pi({\cal A},t)=(1+t)(1+3t)^2$. This is the simplest example of a 
nonfree arrangement where $\pi({\cal A},t)$ factors. Failure of
freeness can be explained using the Addition-Deletion theorem of
\cite{ot}; \lemref{lem:bigmu} provides another explanation.
We combine the results above into 
\end{exm}
\begin{thm}\label{thm:constraints}
Let ${\cal A}$ be an arrangement on $d$ lines with intersection 
lattice $L_{\cal A}$, $M = \max \{ \mu(x)  |  x \in L_2({\cal A}) \}$
 and $\pi({\cal A},t)=(1+t)(1+at)(1+(d-1-a)t)$.
For $$\{ \alpha_1 \le \alpha_2 \le \cdots \le \alpha_m \} \in \mathbb{N}^m
\mbox{ and }\{ \beta_1 \le \beta_2 \le \cdots \le \beta_{m-2} \}\in \mathbb{N}^{m-2},$$
a unique integral solution to the following inequalities implies 
that Terao's conjecture holds for arrangements with intersection 
lattice $L_{\cal A}$.
\begin{enumerate}
\item The bound on the number of generators of $D_0$:
$$2 \le m \le d-1$$
\item Global geometric constraints on the Jacobian ideal:
$$\sum\limits_{i=1}^m \alpha_i - \sum\limits_{j=1}^{m-2} \beta_j = d-1$$
$${d-2 \choose 2} -\sum\limits_{i=1}^m {\alpha_i-1 \choose 2}
+\sum\limits_{j=1}^{m-2} {\beta_j-1 \choose 2} +1=a(d-1-a).$$
\item Regularity constraints:
$$2 \le \alpha_i \le d-2 $$
$$3 \le \beta_j \le d-1.$$
\item Resolution constraints:
$$\beta_1 \ge \alpha_2+1.$$
\item Local geometric constraints:
$$ \exists \mbox{ }\alpha_i \ge M.$$
\end{enumerate}
\end{thm}
\begin{exm}\label{exm:nonfano}
The non-Fano arrangement is an arrangement of seven lines with
$\pi({\cal A},t)=(1+t)(1+3t)^2$:

\begin{center}
\epsfig{file=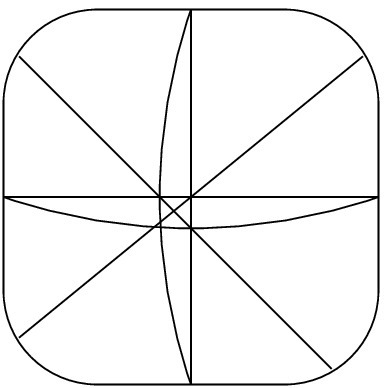,height=1.1in,width=1.1in}
\end{center}
It is the smallest free arrangement which is not supersolvable. 
Applying \thmref{thm:constraints}, we find that there are nineteen 
different numerical possibilities for the free resolution of $D_0$. 
To shorten this list, we employ stability. 
\end{exm}
Write ${\cal A} \equiv {\cal B}$ if ${\cal A}$ and ${\cal B}$ have
isomorphic intersection lattices, and let ${\cal D}_{\cal A}$,
${\cal D}_{\cal B}$ be the bundles associated to arrangements 
${\cal A}$ and ${\cal B}$. Call a split two bundle ${\cal F}$ 
{\it balanced} if 
${\cal F}\simeq {\cal O}^2_{\mathbb{P}^2}(a)$; a rank two bundle which 
splits is semistable iff it is balanced. 
It is obvious that the only free line arrangements which can be semistable  
are those with an odd number of lines, so for the remainder of this section
we'll assume the number of lines is odd. Finally, note that the
normalization of a balanced split two bundle is just 
${\cal O}_{\mathbb{P}^2}$; so both Chern classes of 
the normalized bundle are zero. 
\begin{lem}\label{lem:stab}
Suppose ${\cal D}_{\cal A}$ is balanced, and ${\cal B}\equiv{\cal A}$ 
is a counterexample to Terao's conjecture. Then 
${\cal D}_{\cal B}$ is not semistable. In particular,
there must be a syzygy on the Jacobian ideal of ${\cal B}$ of degree 
$< (d-1)/2$.
\end{lem}
\begin{proof}
Once we've fixed a lattice, we have also fixed the Chern classes. Thus, a
counterexample to
Terao's conjecture, if semistable, would have a jump locus of degree 
$c_2 =0$ by Barth's theorem, so a semistable counterexample would be 
a uniform bundle. But a uniform two bundle on $\mathbb{P}^2$ which does not 
split is ${\cal T}^1(a)$ (\cite{v}), and these possibilities are excluded 
simply by considering the Chern classes.
So we know that a counterexample must be unstable. 
Now from \lemref{lem:norm} we see that $H^0({\cal D}_{{\cal B} norm}(-1)) \ne 0$, 
which implies the result about the syzygies.
\end{proof}
\begin{lem}\label{lem:gbgs}
Any counterexample ${\cal B}$ $(\equiv{\cal A}$ with ${\cal D}_{\cal A}$ balanced$)$ to Terao's 
conjecture must have a syzygy on the Jacobian ideal of
${\cal B}$ of degree $>(d-1)/2$.
\end{lem}
\begin{proof}
If not, then ${\cal D}_{{\cal B} norm}$ would be generated by global sections.
But a globally generated bundle with $c_1 = 0$ is trivial (\cite{oss} p. 53).
\end{proof}
Combining these lemmas yields the following theorem. Notice that since having
syzygies of low degree corresponds to being in special position, this
is consistent with Yuzvinsky's results.
\begin{thm}\label{thm:unstable}
Suppose a free arrangement ${\cal A}$ on $d$ lines has 
${\cal D}_{\cal A}\simeq {\cal O}^2((d-1)/2)$. A counterexample ${\cal B} \equiv {\cal A}$ 
 to Terao's conjecture must have a syzygy of degree $< (d-1)/2$, and also
a syzygy of degree $> (d-1)/2$.
\end{thm}
\begin{exm}\label{thm:nonf}
We return to case where ${\cal A}$ is the non-Fano arrangement. 
Of the nineteen numerically possible free resolutions for $D_0$, 
one corresponds to the case where ${\cal A}$ is
free. \thmref{thm:unstable} allows us to rule out fifteen of the other
possibilities. Thus, there are only three numerical types of resolution 
possible if ${\cal A}$ is not free:
\begin{center}
$ 0 \longrightarrow R(-6) \longrightarrow
R^2(-5)\oplus R(-2) \longrightarrow D_0 \longrightarrow 0.$
\end{center}
\begin{center}
$ 0 \longrightarrow R(-6)\oplus R^3(-3) \longrightarrow
R(-5)\oplus R^3(-4)\oplus R^2(-2) \longrightarrow D_0 \longrightarrow
0.$
\end{center}
\begin{center}
$ 0 \longrightarrow R^2(-6)\oplus R^2(-3) \longrightarrow
R^4(-5) \oplus R^2(-2) \longrightarrow D_0 \longrightarrow 0.$
\end{center}
%
%
%
%
If there are only two quadratic first syzygies, as in the last two
cases, there can be at most one linear second syzygy, so these
resolutions cannot occur. The first possibility actually does occur as the
free resolution for an arrangement on seven lines with Poincar\'e 
polynomial $\pi({\cal A},t)=(1+t)(1+3t)^2$ - it is the free resolution
for the arrangement in Example $5.3$. However, it cannot 
be the free resolution for
the non-Fano arrangement. To see this, start out with arrangement $I$
from last section, and apply \thmref{thm:cmreg} twice. The regularity can
be at most four, whereas this resolution has regularity five. Hence,
any arrangement combinatorially equivalent to the non-Fano arrangement
must be free. Remark: Theorem 4.51 of \cite{ot} and 
a similar sequence of subarrangements can also be used to obtain this
result.
\end{exm}
\noindent{\bf Acknowledgement} I thank two anonymous referees for useful
suggestions. The Macaulay2 software package was used to perform all
computations. 

\bibliographystyle{amsalpha}

\end{document}